\documentclass[a4paper, 11pt]{amsart}
\usepackage{mathpazo}
\usepackage{amsmath, amssymb, amsfonts, amsthm, graphicx, overpic}

\newcommand{\epsi}{\varepsilon}
\newcommand{\norm}[1]{\lVert#1\rVert}
\newcommand{\R}{\mathbf{R}}
\newcommand{\collection}[1]{{\mathcal#1}}
\newcommand{\CK}{\collection{K}}
\newcommand{\abs}[1]{\lvert#1\rvert}

\theoremstyle{plain}
\newtheorem{theorem}{Theorem}
\newtheorem{lemma}{Lemma}
\newtheorem{corollary}{Corollary}
\newtheorem{conjecture}{Conjecture}
\newtheorem{proposition}{Proposition}

\begin{document}

\title[Quantitative illumination and Steiner minimal trees]{Quantitative illumination
of convex bodies and vertex degrees of geometric Steiner minimal trees}
\author{Konrad J. Swanepoel}
\address{Department of Mathematical Sciences,
        University of South Africa, PO Box 392,
        Pretoria 0003, South Africa}
\email{\texttt{swanekj@unisa.ac.za}}

\begin{abstract}
In this note we prove two results on the quantitative illumination parameter $f(d)$ of 
the unit ball of a $d$-dimensional normed space introduced by K.~Bezdek (1992).
The first is that $f(d)=O(2^d d^2\log d)$.
The second involves Steiner minimal trees. 
Let $v(d)$ be the maximum degree of a vertex, and $s(d)$ of a Steiner point, in a 
Steiner minimal tree in a $d$-dimensional normed space, where both maxima are over all norms.
F.~Morgan (1992) conjectured that $s(d)\leq 2^d$, and D.~Cieslik (1990) 
conjectured $v(d)\leq 2(2^d-1)$.
We prove that $s(d)\leq v(d)\leq f(d)$ which, combined with the above estimate of $f(d)$, improves the previously best known upper bound 
$v(d)<3^d$.
\end{abstract}

\maketitle

\section{Introduction}
Let $K$ denote a convex body in the $d$-dimensional real vector space $\R^d$.
Denote its volume by $\mu(K)$ and its translative covering density  by $\vartheta(K)$.
A (positive) \emph{homothet} with ratio $\lambda>0$ of $K$ is any set of the form 
$\lambda K+t$, with $t\in\R^d$.
The \emph{difference body} of $K$ is $K-K$.
According to the Rogers-Shephard inequality \cite{MR19:1073f}, 
$\mu(K-K)/\mu(K)\leq \binom{2d}{d}$.
If $K$ is \emph{centred} (that is, $K=-K$), then of course $\mu(K-K)/\mu(K)=2^d$,
and $K$ defines a norm
\[ \norm{x}_K := \inf\{\lambda >0: \lambda^{-1}x\in K\},\]
which turns $\R^d$ into a normed space.
Let $\CK^d$ denote the class of all $d$-di\-men\-sion\-al convex bodies,
and $\CK^d_o$ the class of all centred $d$-dimensional convex bodies.

\subsection{Quantitative illumination and covering}
A point $p\notin K$ \emph{illuminates} a point $q$ on the boundary of $K$ if the ray
\[ \{\lambda p+(1-\lambda) q : \lambda<0\} \]
intersects the interior of $K$.
A set of points $P\subseteq\R^d\setminus K$ \emph{illuminates} $K$ if each boundary point
of $K$ is illuminated by some point in $P$.
Let $L(K)$ be the smallest size of a set that illuminates $K$.
Also let $L(d):=\max \{L(K):K\in\CK^d\}$, and $L_o(d):=\max \{L(K):K\in\CK^d_o\}$.
Since $L(K)=2^d$ if $K$ is a cube, $L(d)\geq L_o(d)\geq 2^d$.
The well-known illumination problem is to show that $L(d)= 2^d$.
For large $d$ the best known upper bounds are
$L(d)\leq \binom{2d}{d}d(\log d+\log\log d+5)$
and
$L_o(d)\leq 2^dd(\log d+\log\log d+5)$,
due to Rogers \cite[p.~284]{MR27:4134}; see also \cite{MR98i:52026}.

There are other equivalent formulations of this illumination problem.
For example, let $L'(K)$ be the smallest number of positive homothets of $K$,
with each homothety ratio less than $1$, whose union contains $K$.
Then $L(K)=L'(K)$.
See \cite{MR2000b:52006} for a survey on this problem and its history.

We consider quantitative versions of the above two formulations of 
the illumination problem.
The first was introduced by K.~Bezdek \cite{Bezdek92}.
For $K\in\CK^d_o$ let
\[B(K):=\inf\left\{\sum_i\norm{p_i}_K : \{p_i\} \text{ illuminates } K\right\}.\]
This ensures that far-away light sources are penalised.
Let \[B(d):=\sup \{B(K): K\in\CK^d_o\}.\]
Bezdek asked for the value of $B(d)$, and in particular, if $B(d)$ is finite
for $d\geq 3$.
He showed that $B(2)=6$; the regular hexagon giving equality.
Note that $B(K)\geq L(K)$, hence $B(d)\geq L_o(d)\geq 2^d$.
It is also easily seen that $B(K)=2^d$ if $K$ is a $d$-cube, and $B(K)=2d$ if
$K$ is a $d$-cross polytope.

We introduce the following quantitative covering parameter for $K\in\CK^d$:
\[ C(K):= \inf\left\{\sum_i(1-\lambda_i)^{-1}: 
K\subseteq\bigcup_i(\lambda_i K+t_i), 0<\lambda_i<1, t_i\in\R^d\right\}.\]
In this way homothets almost as large as $K$ are penalised.
\begin{proposition}\label{propa}
For any $K\in\CK^d_o$ we have $B(K)\leq 2C(K)$.
\end{proposition}
Let \[C(d):=\sup \{C(K):K\in\CK^d\},\] and \[C_o(d):=\sup \{C(K):K\in\CK^d_o\}.\]
Hence $C(d)\geq C_o(d)\geq B(d)/2$.
It is easy to see that $C(K)=2^{d+1}$ if $K$ is a $d$-cube, hence
$C(d)\geq C_o(d)\geq2^{d+1}$.
As before, it is not clear whether $C(d)$ is finite.
Levi \cite{MR16:163a} showed that any planar convex body can be covered with $7$ homothets,
each with homothety ratio $1/2$;
 hence $C(2)\leq 14$.
Lassak's result \cite{Las86} that any planar convex body can be covered with $4$ homothets,
each with ratio $1/\sqrt{2}$, improves this to $C(2)\leq 8+4\sqrt{2}$.
Lassak \cite{MR99f:52023} also showed that any convex body in $\R^3$ can be covered with 
$28$ homothets, each with ratio $7/8$;
 hence $C(3)\leq 224$.
We show that a result of Rogers and Zong \cite{MR98i:52026} implies the following
upper bound.
\begin{theorem}\label{prop1}
For any $d$-dimensional convex body $K$ we have
\[C(K) < e(d+1)\frac{\mu(K-K)}{\mu(K)}\vartheta(K).\]
\end{theorem}
Using Rogers' estimate \cite{MR19:877c} $\vartheta(K)\leq d(\log d+\log\log d+5)$
for $d\geq 2$ and the Rogers-Shephard inequality one finds
\[C(d) < \binom{2d}{d}e(d+1)d(\log d+\log\log d+5)=O(4^d d^{3/2}\log d),\] and
\[ B(d) \leq 2C_o(d) < 2^{d+1}e(d+1)d(\log d+\log\log d+5)=O(2^d d^2\log d).\]
Perhaps $C(d)=O(2^d)$.

\subsection{Steiner minimal trees}
Given a finite set of points $V$ in $\R^d$, a \emph{Steiner tree} $T$ of $V$ is any tree
in $\R^d$ whose vertex set contains $V$, and whose edges are straight-line segments in $\R^d$.
The vertices of $T$ not in $V$ are called \emph{Steiner points}.
(Usually Steiner points are required to have degree at least $3$, but this is 
unnecessary here.)
The \emph{$K$-length} of a Steiner tree is the total length in $\norm{\cdot}_K$ of the
edges of the tree, where $K$ is a centred convex body.
It is easily seen \cite{Cockayne} that any given point set has a Steiner tree of smallest $K$-length,
called a \emph{$K$-Steiner minimal tree} ($K$-SMT).

Steiner minimal trees have been studied mostly in the Euclidean plane and the rectilinear
plane ($K$ a parallelogram) \cite{MR94a:05051}.
Other normed planes have also been considered; see \cite[\S 3.1]{MR2003m:90136} for
further references.
Steiner minimal trees in normed spaces of higher dimension have been investigated by Cieslik
\cite{MR99i:05062} and Morgan \cite{MR93h:53012} among others.

Let $v(K)$ be the maximum possible degree of a vertex in a $K$-SMT, and $s(K)$ the
maximum possible degree of a Steiner point in a $K$-SMT.
Clearly $s(K)\leq v(K)$.
The following table gives some examples of known values of $s(K)$ and $v(K)$.
See \cite{MR2000g:05054, MR2001f:05049, MR2001m:05084} for further examples.

\medskip
\begin{center}
\begin{tabular}{|c|c|c|} \hline
$K$ & $s(K)$ & $v(K)$ \\ \hline
Euclidean $d$-ball & $3$ & $3$ \\
$d$-cube & $2^d$ & $2^d$ \\
$d$-cross polytope & $2d$ & $2d$ \\
regular hexagon & $4$ & $6$ \\ \hline
\end{tabular}
\end{center}

\medskip
Let $s(d):=\max \{s(K):K\in\CK^d_o\}$, and $v(d):=\max \{v(K):K\in\CK^d_o\}$.
Then $2^d\leq s(d)\leq v(d)$.
The following two conjectures have been made:
\begin{conjecture}[Cieslik \cite{MR92a:05039}, {\cite[ch.~4]{MR99i:05062}}]
$v(d)\leq 2(2^d-1)$ for all $d\geq 2$. \label{conj1}
\end{conjecture}

\begin{conjecture}[Morgan \cite{MR93h:53012}, {\cite[ch.~10]{MR98i:53001}}]
$s(d)\leq 2^d$ for all $d\geq 2$. \label{conj2}
\end{conjecture}

Cieslik \cite{MR92a:05039} has shown that $v(K)\leq H(K)$ where $H(K)$ is the
translative kissing number of $K$.
See \cite{MR98k:52048} for a survey and for references to the following upper bounds on $H(K)$.
Since $H(K)\leq 3^d-1$ with equality only for (affine images of) the $d$-cube, 
it follows that $v(d)\leq 3^d-2$.
Since for planar $K$ we have $H(K)\leq 6$ if $K$ is not a parallelogram,
we obtain $v(2)= 6$ \cite{MR92a:05039}; thus Conjecture~\ref{conj1} is true for $d=2$.
Conjecture~\ref{conj2} is also true for $d=2$ \cite{MR2001f:05049}.
The two-dimensional methods are very special and offer no hope for generalisation to
higher dimensions.
We find upper bounds within a factor of $O(d^2\log d)$ from the conjectured values,
using the following relationship with Bezdek's illumination parameter.
\begin{theorem}\label{thm2}
For any $K\in\CK^d_o$ we have $v(K)\leq B(K)$.
\end{theorem}
Note that equality holds, for example, if $K$ is a regular hexagon, a $d$-cube or
a $d$-cross polytope, but not if $K$ is a $d$-ball.
\begin{corollary}
For any $K\in\CK^d_o$ we have $s(K)\leq v(K) < 2^{d+1} e(d+1)\vartheta(K)$.
\end{corollary}
\begin{corollary}
$s(d)\leq v(d) =O(2^d d^2\log d)$.
\end{corollary}

\section{Proofs}
\begin{proof}[Proof of Proposition~\ref{propa}]
Let $\{\lambda_i K+t_i\}$ be a finite covering of $K$, with $0<\lambda_i<1$ for all $i$.
Let $\epsi>0$  be sufficiently small such that all $\lambda_i+\epsi<1$.
If a boundary point $q$ of $K$ is covered by $\lambda_i K+t_i$, then
$1-\lambda_i \leq \norm{t_i}_K \leq 1+\lambda_i<2$, and the centre of the homothety mapping
$K$ to $(\lambda_i+\epsi) K+t_i$, namely $p_i:=(1-\lambda_i-\epsi)^{-1}t_i$, is outside $K$
and illuminates $q$.
Therefore, the set $\{p_i\}$ illuminates $K$, and
$\sum_i\norm{p_i}_K<\sum_i 2/(1-\lambda_i-\epsi)$.
Since $\epsi>0$ can be made arbitrarily small,
$\sum_i\norm{p_i}_K\leq2\sum_i(1-\lambda_i)^{-1}$.
\end{proof}

\begin{proof}[Proof of Theorem~\ref{prop1}]
It is known \cite{MR98i:52026} that for any $0<\lambda<1$ there exists a covering of
$K$ by homothets $\{\lambda K+t_i : i=1,\dots,N\}$, with
\[ N \leq \frac{\mu(K-\lambda K)}{\mu(\lambda K)}\vartheta(K)<\lambda^{-d}\frac{\mu(K-K)}{\mu(K)}\vartheta(K).\]
Choosing $\lambda=d/(d+1)$ we obtain
\[ \sum_{i=1}^N(1-\lambda)^{-1}
< (d+1)\left(1+\frac{1}{d}\right)^d\frac{\mu(K-K)}{\mu(K)}\vartheta(K)
< (d+1)e \frac{\mu(K-K)}{\mu(K)}\vartheta(K).\]
\end{proof}

\begin{lemma}\label{lemma1}
If $p$ illuminates the boundary point $u$ of $K\in\CK^d_o$, then for all sufficiently small
$\epsi>0$, $\norm{u-\epsi p}_K<1-\epsi$.
\end{lemma}
\begin{proof}
The lemma is trivial if $p=\lambda u$ for some $\lambda$.
Therefore, assume that $p$ and $u$ are linearly independent and consider the two-dimensional
subspace spanned by them (Figure~\ref{fig1}).
\begin{figure}[h]
\begin{overpic}[scale=0.6, clip=true, bb=0 120 383 253]{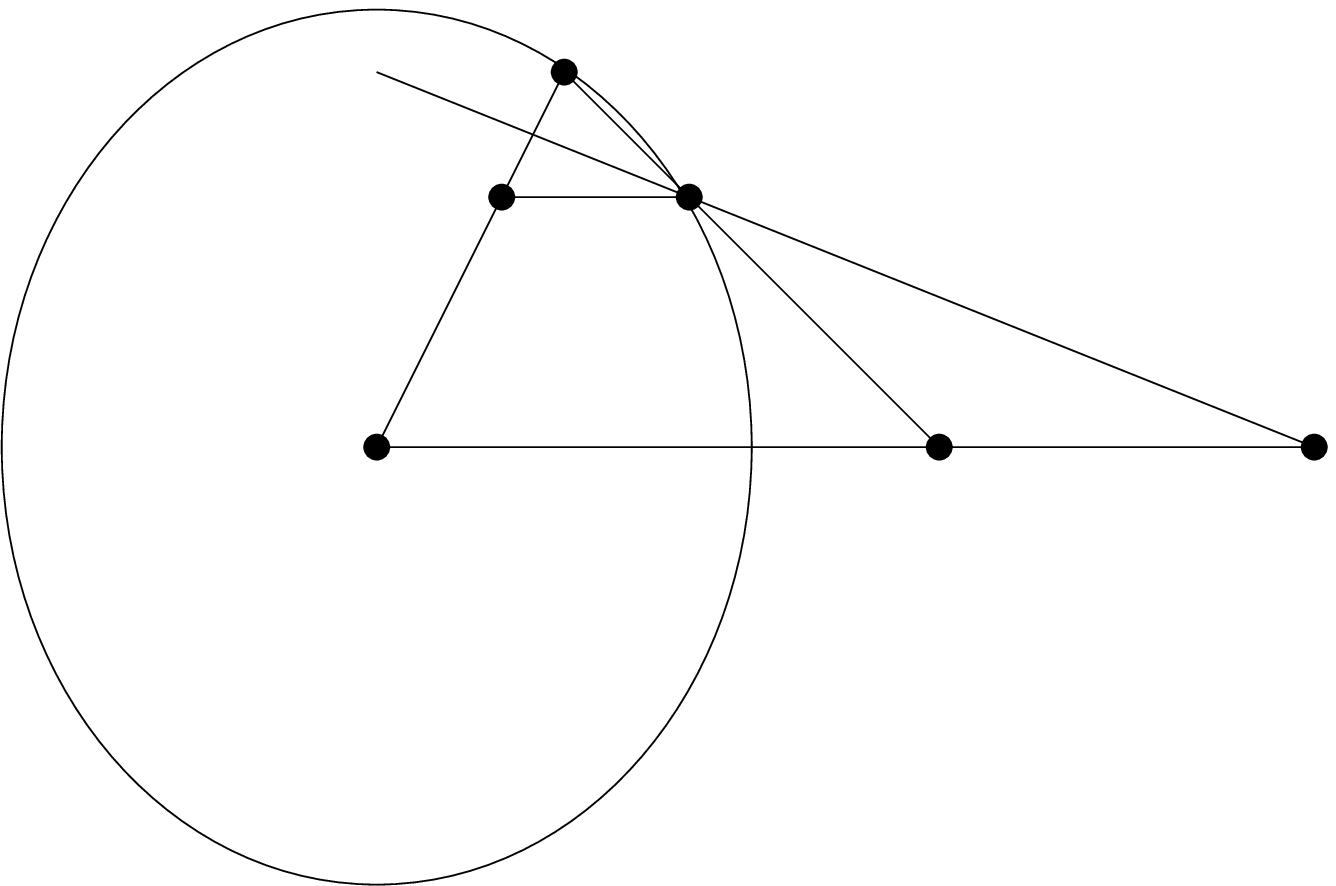}
\put(44,32){$v$}
\put(53,22){$u$}
\put(22,20){$u-\epsi p$}
\put(39,15){$\norm{\epsi p}_K$}
\put(70,15){$\ell$}
\put(32,3){$o$}
\put(71,5){$p'$}
\put(98,5){$p$}
\end{overpic}
\caption{}\label{fig1}
\end{figure}
Since $p$ illuminates $u$, we may choose $\epsi_0>0$ such that the line through $o$ and
$u-\epsi_0 p$ intersects the line $\ell$ through $u$ and $p$ in the interior of $K$.
Then clearly for all $\epsi>0$ with $\epsi<\epsi_0$ the line through $o$ and $u-\epsi p$
still intersects $\ell$ in the interior of $K$.
Let $v=(\norm{u-\epsi p}_K)^{-1}(u-\epsi p)$.
Then the lines $vu$ and $op$ intersect in $p'$, say, with $\norm{p'}_K<\norm{p}_K$.
Using similar triangles, $\norm{u-\epsi p}_K=1-\norm{\epsi p}_K/\norm{p'}_K < 1-\epsi$.
\end{proof}

\begin{proof}[Proof of Theorem~\ref{thm2}]
Consider a vertex of a $K$-SMT of degree $v(K)$.
By translating we may assume that the vertex is the origin $o$.
By scaling we may also assume that each edge emanating from $o$ has $K$-length at least $1$.
Let these edges be $ov_i$, with $\norm{v_i}_K\geq 1$.
Let $u_i=\norm{v_i}_K^{-1}v_i$.
Then the star $T$ joining $o$ to each $u_i$ is a $K$-SMT of $\{o,u_1,u_2,\dots,u_{v(K)}\}$
(otherwise we would be able to shorten the original tree).

Let $\{p_1,\dots,p_k\}$ illuminate $K$.
For each $j=1,\dots,k$, let \[U_j=\{u_i: p_j \text{ illuminates } u_i\}.\]
Then $\{u_i\}=\bigcup_j U_j$.
We estimate the number of points $\abs{U_j}$ in each $U_j$.
By Lemma~\ref{lemma1} we may find $\epsi>0$ such that $\norm{u_i-\epsi p_j}_K<1-\epsi$
for all $i$.
Consider the tree $T'$ obtained from the star $T$ by replacing, for each $u_i\in U_j$,
the edge from $o$ to $u_i$ by the edge from $\epsi p_j$ to $u_i$, and joining the
Steiner point $\epsi p_j$ to $o$.
Then $T'$ is not shorter than $T$.
This implies that
\begin{eqnarray*}
\abs{U_j}=\sum_{u_i\in U_j}\norm{u_i}_K
&\leq & \norm{\epsi p_j}_K + \sum_{u_i\in U_j}\norm{u_i-\epsi p_j}_K\\
&<& \epsi\norm{p_j}_K + (1-\epsi)|U_j|,
\end{eqnarray*}
and $|U_j|<\norm{p_j}_K$.
Hence $v(K)\leq\sum_{j=1}^k|U_j| <\sum_{j=1}^k\norm{p_j}_K$.
Taking the infimum over all sets $\{p_i\}$ that illuminate $K$, we obtain that
$v(K)\leq B(K)$.
\end{proof}

\renewcommand{\MR}[1]{\relax}

\end{document}